\documentclass[12pt]{amsart}
\usepackage{amssymb, epsfig,amsmath}
\usepackage[mathscr]{euscript}
\usepackage{amscd}
\usepackage{stmaryrd}

\newcommand{\bba}{\mathbb{A}}
\newcommand{\bbc}{\mathbb{C}}

\newcommand{\bbq}{\mathbb{Q}}
\newcommand{\bbr}{\mathbb{R}}
\newcommand{\bbz}{\mathbb{Z}}
\newcommand{\bbf}{\mathbb{F}}
\newcommand{\bbn}{\mathbb{N}}
\newcommand{\bbh}{\mathbb{H}}

\newcommand{\Fn}{\mathfrak{n}}
\newcommand{\Fp}{\mathfrak{p}}
\newcommand{\Fq}{\mathfrak{q}}

\newcommand{\CO}{\mathcal{O}}
\newcommand{\CN}{\mathcal{N}}

\newcommand{\BCinf}{\bbc_{\infty}}

\DeclareMathOperator{\Gal}{Gal}
\DeclareMathOperator{\End}{End}
\DeclareMathOperator{\Pic}{Pic}
\DeclareMathOperator{\GL}{GL}

\newcommand{\ks}{\overline{k}^{\sigma}}
\DeclareMathOperator{\tor}{tors}

\DeclareMathOperator{\id}{id}

\newtheorem{thm}{Theorem}[section]

\newtheorem{lem}[thm]{Lemma}
\newtheorem{prop}[thm]{Proposition}
\newtheorem{conj}[thm]{Conjecture}

\newtheorem{defnn}[thm]{Definition}

\newtheorem{remarkk}[thm]{Remark}

\newtheorem{examplee}[thm]{Example}

\title[The rank of elliptic curves over global fields]
{Heegner points and the rank of elliptic curves over large
extensions of global fields}
\author{Florian Breuer and Bo-Hae Im}
\date{January 31, 2006}
\address{Department of Mathematical Sciences, University of Stellenbosch, Stellenbosch 7600, South Africa}
\email{fbreuer@sun.ac.za}
\address{Department of Mathematics, University of Utah, Salt Lake
City, Utah 84112, USA} \email{im@math.utah.edu}
\subjclass[2000]{Primary 11G05}

\begin{document}

\begin{abstract} Let  $k$ be a global field, $\overline{k}$ a separable
closure of $k$, and $G_k$ the absolute Galois group
$\Gal(\overline{k}/k)$ of $\overline{k}$ over $k$. For every
$\sigma\in G_k$, let $\ks$ be the fixed subfield  of $\overline{k}$
under $\sigma$. Let $E/k$ be an elliptic curve over $k$. We show
that for each $\sigma\in G_k$, the Mordell-Weil group $E(\ks)$ has
infinite rank in the following two cases. Firstly when $k$ is a
global function field of odd characteristic and $E$ is parametrized
by a Drinfeld modular curve, and secondly when $k$ is a totally real
number field and $E/k$ is parametrized by a Shimura curve. In both
cases our approach uses the non-triviality of a sequence of Heegner
points on $E$ defined over ring class fields.
%
\end{abstract}

\maketitle

\vspace{1 cm}
\section{Introduction}
 This paper is
motivated by the following conjecture of M. Larsen (see \cite{larsen} for context):

\vspace{.5 cm}

 \noindent{\bf Conjecture.}\hspace{.1 cm} Let $E/k$ be an elliptic curve over
 a finitely generated infinite field $k$. Then, for every
$\sigma\in G_k := \Gal(\overline{k}/k)$, the Mordell-Weil group
$E(\overline{k}^{\sigma})$ of $E$ over
 $\overline{k}^{\sigma}=\{x\in\overline{k}\mid \sigma(x)=x\}$
  has infinite rank.
\vspace{.5 cm}

\noindent Larsen \cite{larsen} has shown that the conjecture is true for $\sigma$ in
a suitable open subset of $G_k$.
In \cite{im}, \cite{im2} and \cite{im3}, the conjecture is proved for the following
cases over number fields:
%
 \begin{enumerate}
 \item If $E/k$ has a $k$-rational point $P$ such that $2P\neq O$ and $3P\neq O$, or
 \item if the 2-torsion points of $E/k$ are $k$-rational, or
 \item if $k=\bbq$ without any hypothesis on rational points,
   \end{enumerate}
then, for every automorphism $\sigma\in G_k$, the rank of  the
Mordell-Weil group
 $E(\overline{k}^{\sigma})$
is infinite.

In fact, under the first assumption above
it is shown in \cite{im} that for each
$\sigma\in\Gal(\overline{k}/k)$, the rank of $E$ over the maximal
Galois extension of $k$ in $\overline{k}^\sigma$ is infinite. In
general, the maximal Galois extension of $k$ in
$\overline{k}^\sigma$ is smaller than $\overline{k}^\sigma$.
Moreover,  under the second assumption above that the 2-torsion
points of $E/k$ are $k$-rational, it was shown in \cite{im2} that
$E$ has infinite rank over the maximal abelian extension of $k$ in
$\overline{k}^\sigma$, which is also much smaller.

In this paper, we prove that Larsen's conjecture is true for certain
modular elliptic curves over global fields. When $k$ is a totally
real number field, We say $E/k$ is modular if $E/k$ is parametrized
by a suitable Shimura curve. Every elliptic curve over $k=\bbq$ is
parametrized by a modular curve $X_0(N)$, where $N$ is the conductor
of $E/\bbq$, by \cite{Wi}, \cite{TW} and \cite{BCDT}. $X_0(N)$ is a
particularly simple Shimura curve, so $E/\bbq$ is also modular in
the above sense. When $k$ is a global function field and $E/k$ has
split multiplicative reduction at a place $\infty$, then the
conductor of $E/k$ may be written as $\Fn\cdot\infty$, where $\Fn$
is an ideal in the Dedekind ring $A:=\{x\in k\;|\; \text{$x$ is
regular away from $\infty$}\}$, and $E/k$ is parametrized by the
Drinfeld modular curve $X_0(\Fn)$, see \cite{GR}.

%

Our approach is the following. Let $E/k$ be a modular elliptic curve, then for a given
automorphism $\sigma\in G_k$, we produce an
infinite sequence $K_m$ of imaginary quadratic extensions of $k$
(if $k$ is a function field, then $K_m/k$ is {\em imaginary} if $\infty$ does not split in $K_m$) such that
\begin{enumerate}
\item if $\sigma|_{K_m}=\id_{K_m}$ for all $m$, then we may show
by elementary arguments that the rank of $E(\ks)$ is infinite
and
\item if $\sigma|_{K_m}\neq \id_{K_m}$ for some $m$, then $(E,
K_m)$ satisfies the \emph{Heegner hypothesis} (the definition is
given below). This allows us to construct a suitable sequence of Heegner points on
$E$ defined over a tower of ring class fields of $K_m$ over which the rank of
$E$ is unbounded. Then we use the dihedral structure of these ring class fields to show
that the rank of $E(\ks)$ is infinite.
\end{enumerate}

Note that for number fields we simplify the proof of the results of \cite{im3}, using an
argument from \cite{breuer2}.

If $k$ is a function field and $E/k$ has non-constant $j$-invariant, then there exists a finite extension
$L/k$ such that $E/L$ has split multiplicative at some place of $L$,
and hence our results apply to $E/L$, but this case is already covered by \cite[Theorem 5]{larsen}.

\vspace{.5 cm} \noindent {\bf Acknowledgements. } We would like to
thank Henri Darmon for suggesting the case of totally real fields of
this problem and for his encouragement. Also we would like to thank
Michael Larsen for for his interest in our result and his
encouragement.

\vspace{.5 cm}

\section{Modular elliptic curves over global function fields}

 Throughout this section, we let $k$ be a global function field with
field of constants $\bbf_q$, where $q$ is a power of the odd
prime $p$, and we let $A$ be the ring of functions in $k$ regular outside the place $\infty$.
Let $E/k$ be an elliptic curve over $k$ with split multiplicative
reduction at $\infty$, then the conductor of $E/k$ is $\Fn\cdot\infty$ for an ideal $\Fn\subset A$.

Denote by $X_0(\Fn)$ the Drinfeld modular curve parametrizing pairs of rank-2 Drinfeld $A$-modules
linked by cyclic $\Fn$-isogenies. Then we have a morphism defined over $k$ (see \cite{GR}):
\[
\pi: X_0(\Fn) \longrightarrow E.
\]

Let $K/k$ be a quadratic {\em imaginary} extension, i.e. $\infty$ does not split in $K/k$.
Suppose that all prime divisors of $\Fn$ split in $K/k$, we say that $(E,K)$ satisfies the
{\bf Heegner hypothesis}.

We fix a prime $\Fp$ of $A$ not dividing $\Fn$.

Denote by $\CO_K$ the integral closure of $A$ in $K$, then thanks to the Heegner hypothesis, there exists
an ideal $\CN\subset\CO_K$ such that $\CO_K/\CN\cong A/\Fn$. For every non-negative integer $n$
we denote by
$\CO_n := A + \Fp^n\CO_K$ the order in $\CO_K$ of conductor $\Fp^n$, and set $\CN_n:=\CN\cap\CO_n$. We have
$\CO_n/\CN_n\cong A/\Fn$ for all $n\geq 0$. Denote by $\BCinf = \hat{\bar{k}}_{\infty}$ the completion
of an algebraic
closure of the completion of $k$ at $\infty$, a field both algebraically closed and complete,
which plays the role
of the complex numbers in characteristic $p$. Then $\CO_K$ and $\CN_n^{-1}$ are rank-2 lattices in $\BCinf$,
hence define
a pair of Drinfeld modules $(\Phi^{\CO_K}, \Phi^{\CN_n^{-1}})$ linked by a cyclic $\Fn$-isogeny.
The pair thus
defines a point $x_n$ on $X_0(\Fn)$, which is defined over the ring class field $K[\Fp^n]$
of the order $\CO_n$.
Its image $y_n=\pi(x_n)\in E(K[\Fp^n])$ is called a {\em Heegner point} on $E$.

Recall that $\Gal(K[\Fp^n]/K)\cong\Pic(\CO_n)$.
We now show that the rank of $E(K[\Fp^n])$ is unbounded as $n\rightarrow\infty$.
This follows from \cite[Theorem 4.4]{breuer3}, but notice that in that paper and in \cite{breuer}
the Heegner points $y_n$
are constructed differently, involving a trace, and consequently the proof is more difficult.
For our purposes the following result for our Heegner points is sufficient, and much easier.

\begin{prop}\label{infrank}
Let $K[\Fp^\infty] = \bigcup_{n\geq 0}K[\Fp^n]$. Let $I\subset\bbn$ be an infinite set.
Then the subgroup of $E(K[\Fp^\infty])$ generated by $\{y_n \;|\; n\in I\}$ has finite torsion and
infinite rank.
\end{prop}

\begin{proof}
From \cite[Lemma 2.2]{breuer} follows that $E(K[\Fp^\infty])$ has
finite torsion, so by \cite[Lemma~2.5]{im3}, it remains to show that
the subgroup $H\subset E(K[\Fp^\infty])$ generated by the $y_n$'s is
not finitely generated. For this we use the same argument as in
\cite{breuer2}.

Suppose that $H$ is finitely generated. Then $H\subset E(L)$ for some finite separable extension $L/k$,
which we may
extend to include $K$. Denote by $G_L=\Gal(\bar{L}/L)$ the absolute Galois group of $L$. Then $G_L$
acts on the
fibers $\pi^{-1}(y_n)$, and the $G_L$-orbit of $x_n$ is bounded: $\# G_L\cdot x_n \leq \deg(\pi)$.

On the other hand,
\begin{eqnarray*}
\# G_L\cdot x_n & \geq & \#\Pic(\CO_n)/[L:K] \\
                & \geq & \frac{\#\Pic(\CO_K)}{[L:K][\CO_K^\times : \CO_n^\times]}|\Fp^n|(1-|\Fp|^{-1}) \\
                & \geq & \frac{\#\Pic(\CO_K)}{[L:K](q+1)}|\Fp|^n(1-|\Fp|^{-1})
\end{eqnarray*}
by equation (2.5) of \cite{breuer3}. This is unbounded as $n$ gets large, which is a contradiction.
\end{proof}

We point out that everything up to now holds almost verbatim if we
replace the function field $k$ with a number field, provided that
$E/k$ is parametrized by the (classical) modular curve $X_0(N)$. In
particular, this provides a simplified proof of the conclusion
of \cite[Proposition 2.7]{im3}.

We will need the following result later.

\begin{prop}\label{group} For positive integers $n$ and $m$ such that $n>m$, the Galois group
$\Gal(K[\mathfrak{p}^n]/K[\mathfrak{p}^m])$ of the ring class fields
$K[\mathfrak{p}^n]$ over $K[\mathfrak{p}^m]$ is a finite
direct sum of cyclic groups of order $p$, where $p$ is the
characteristic of $k$, and the Galois group
$\Gal(K[\mathfrak{p}^n]/k)$ over $k$ is generalized dihedral as in
the number field case.
\end{prop}
\begin{proof} See \cite[Sec.~2.3 and 2.5, Proposition~2.5.7]{brown}.
\end{proof}

%
%

Finally, the following proposition proves the infinite rank of $E$
over $K_{ab}^\sigma$, which will allow us to complete one of our main results,
Theorem~\ref{main}, which is introduced in the next section.

\begin{prop}\label{prop:sigma} Let $\sigma\in G_k$ and $K$ a quadratic imaginary
extension of $k$ such that $\sigma|_K\neq id_K$ and $\infty$ does
not split in $K$. Let $E/k$ be an elliptic curve over $k$ with split
multiplicative reduction at $\infty$. Suppose all primes dividing
$\mathfrak{n}$ split in $K$, where $\mathfrak{n}\cdot\infty$ is the
conductor of $E/k$. Let $\mathfrak{p}\subset A$ be a prime not
dividing $\mathfrak{n}$. Then, the rank of the Mordell-Weil group
$E((K[\mathfrak{p}^n])^\sigma)$ over the fixed subfield of
$K[\mathfrak{p}^n]$ under $\sigma$ is unbounded, as $n$ goes to
$\infty$. In particular, the rank of $E(K_{ab}^\sigma)$ is infinite,
where $K_{ab}$ is the maximal abelian extension of $K$.
\end{prop}

\begin{proof}
%
For the given $\sigma\in G_k$, since $\sigma|_K\neq id_K$, the
restriction of $\sigma$ to each ring class field $K[\mathfrak{p}^n]$
of conductor $n$ can be lifted as an involution of
$K[\mathfrak{p}^n]$. Let $\sigma_n=\sigma|_{K[\mathfrak{p}^n]}$ be
the restriction of $\sigma$ to $K[\mathfrak{p}^n]$. Then, since the
Galois group of  each ring class field $K[\mathfrak{p}^n]$ over $k$
has a generalized dihedral group structure by
Proposition~\ref{group}, $\sigma_n$ acts on $K[\mathfrak{p}^n]$ by
an involution such that
$$(*)\hspace{4 cm} \text{for any }\tau\in \Gal(K[\mathfrak{p}^n]/k),
~~~\sigma_n\tau\sigma_n=\tau^{-1}.\hspace{5 cm}$$

 Now, we prove that the rank of
$E(K[\mathfrak{p}^n]^{\sigma})$ is unbounded as $n$ goes to
infinity.
 Suppose not. Then since the restriction $\sigma_n$ of $\sigma$
 acts by an involution on each
 ring class field $K[\mathfrak{p}^n]$, and
 by Proposition~\ref{infrank} the rank of $E(K[\mathfrak{p}^n])$ is
unbounded as $n$ goes to infinity,
  there exists a fixed integer $n_0$ such that
$\sigma$ acts by $-1$ on any nontrivial quotient
$M_n:=E(K[\mathfrak{p}^n])/E(K[\mathfrak{p}^{n_0}])$, for all $n>
n_0$. Since $K[\mathfrak{p}^{n_0}]$ is Galois over $k$,
$\Gal(K[\mathfrak{p}^{n}]/k)$ acts on $E(K[\mathfrak{p}^{n_0}])$ by
restriction to $K[\mathfrak{p}^{n_0}]$, so it does act on $M_n$.

Let $$\rho : \Gal(K[\mathfrak{p}^{n}]/k) \rightarrow
\mbox{Aut}(M_n)$$ be the representation of
$\Gal(K[\mathfrak{p}^{n}]/k)$. Then, by the hypothesis, $\sigma_n$
acts by $-1$ on $M_n$. Hence, $\rho(\sigma_n)=-id$ on $M_n$.
Therefore,
$$(**)\hspace{1.5cm}\rho(\tau^2)=\rho(\tau)\rho(\tau)=(-id)\rho(\tau)(-id)\rho(\tau)
=\rho(\sigma_n\tau\sigma_n\tau)
\stackrel{by (*)}{=} \rho(1)=id.\hspace{2cm}$$ In particular, for
$\tau\in \Gal(K[\mathfrak{p}^n]/K[\mathfrak{p}^{n_0}])$, the
restriction of $\rho$ to the subgroup $\langle \tau^2\rangle$ of
$\Gal(K[\mathfrak{p}^{n}]/k)$
 generated by the element $\tau^2$ for each
 $\tau\in \Gal(K[\mathfrak{p}^n]/K[\mathfrak{p}^{n_0}])$
 is a
trivial representation of $M_n$. By Proposition~\ref{group},
$\Gal(K[\mathfrak{p}^n]/K[\mathfrak{p}^{n_0}])$ is a $p$-group with
the odd prime $p$ which is the characteristic of $k$. So the order
of every nontrivial element in
$\Gal(K[\mathfrak{p}^n]/K[\mathfrak{p}^{n_0}])$ is the odd prime
$p$. So $\langle \tau^2\rangle=\langle \tau\rangle$. Therefore,
$(**)$ implies that $\rho|_{\langle \tau\rangle}=id$. So we have
that
$$M_n=M_n^{\langle \tau\rangle}, \text{ for all }\tau\in
\Gal(K[\mathfrak{p}^n]/K[\mathfrak{p}^{n_0}]).$$ Since this is true
for all $\tau\in \Gal(K[\mathfrak{p}^n]/K[\mathfrak{p}^{n_0}])$,
this implies that
$$M_n=M_n^{\Gal(K[\mathfrak{p}^n]/K[\mathfrak{p}^{n_0}])}.$$ So
$$E(K[\mathfrak{p}^n]) = E(K[\mathfrak{p}^n]^{\Gal(K[\mathfrak{p}^n]/K[\mathfrak{p}^{n_0}])})
=E(K[\mathfrak{p}^{n_0}]), ~~~\mbox{~~for all }
n > n_0,$$  which is a contradiction to Proposition~\ref{infrank}.

 Therefore,
the rank of $E(K[\mathfrak{p}^n]^{\sigma})$ is unbounded, as
$n\rightarrow\infty$. Since all ring class fields
$K[\mathfrak{p}^n]$ are abelian over $K$, this implies that the rank
of $E((K_{ab})^{\sigma})$ is infinite.
\end{proof}

\section{Infinite rank in odd characteristic $p$}

In this section, we prove our main theorem for the global function
field case. As a generalization of \cite[Lemma]{sil} (which is made
in the number field setting) for global function fields, we will
need the following lemma to prove the linear independence of
algebraic points defined over extensions of bounded degree.

\begin{lem}\label{degreebd} Let $k$ be a global
function field with characteristic $p$ and $E/k$ an elliptic curve
over $k$. Then for any integer $d>1$, the set
$$\bigcup\limits_{[L:k]\leq d} E(L)_{\tor} \text{ is finite.}$$
\end{lem}
\begin{proof} Choose a prime ideal $\pi\subset k$ such that $E$ has a good
reduction at $\pi$. Let $L$ be an extension of degree $\leq d$ over
$k$. Extend $\pi$ to $L$ and denote  its residue field by
$\widetilde{L}$. Then, since $[L:k]\leq d$, $\widetilde{L}$ is
contained in the unique extension $\widetilde{k}_{d!}$ of degree
$d!$ over the residue field $\widetilde{k}$ of $k$ at $\pi$. Note
that for the global function field, a residue field has
characteristic $p$ which is the characteristic of $k$. So we
consider two cases: $n$-torsion points where $p\nmid n$ and
$p^i$-torsion points.

First, for any integer $n$ not divisible by $p$, the $n$-torsion
subgroup $E(L)[n]$ injects into
$\widetilde{E}(\widetilde{k}_{d!})[n]\subset
\widetilde{E}(\widetilde{k}_{d!})$ which doesn't depend on $L$. So
there exists a constant $c_1$ such that  $|E(L)[n]| \leq
|\widetilde{E}(\widetilde{k}_{d!})|=c_1$, for any extension $L$ of
degree $\leq d$ over $k$.

Secondly, for the prime $p$, by \cite[Lemma~1.1]{jac}, we can choose
a prime $\pi\subset k$ such that $E$ has a good reduction at $\pi$
and it induces an isomorphism between $E(\overline{K})[p^i]$ and
$\widetilde{E}\left(\overline{\widetilde{k}}\right)[p^i]$ for all
$i\geq 1$, where $\overline{K}$ is an algebraic closure of $K$. So
again by the same argument in the above, we can show that there
exists a constant $c_2$ such that $|E(L)[p^i]|\leq c_2$ for all $i$
and for any $L$ of degree $\leq d$ over $k$, where $c_2$ does not
depend on $L$ but depends on the degree $d$. Therefore, $
|E(L)_{\tor}| \leq c_1c_2$, for any extension $L$ with $[L:k]\leq
d$. This completes the proof.
\end{proof}

\begin{thm}\label{main} Let $E/k$ be an elliptic curve
over a global function field with the field of constants
$\mathbb{F}_q$, where $q$ is a power of an odd prime $p$. Suppose
$E$ has split multiplicative reduction
at $\infty$. Then, for every automorphism $\sigma\in G_k$, the rank
of $E(\ks)$ is infinite.
\end{thm}

\begin{proof} Let $A$ be the ring of
functions in $k$ regular outside $\infty$. Let
$\mathfrak{n}\cdot\infty$ be the conductor of $E/k$, where
$\mathfrak{n}$ is an ideal in $A$.

Since the characteristic $p$ of $k$ is not 2, we can write a
Weierstrass equation of $E/k$ in the form $y^2=x^3+ax^2+bx+c$. By
a change of variables, we may assume that $a, b$ and $c$ are in
$A$.

Let $\mathfrak{p}_1,\ldots, \mathfrak{p}_m$ be all distinct primes
dividing the ideal $\mathfrak{n}$ in $A$. Fix $M:=p_1\cdots p_m$,
where  $p_i$ is a non-zero prime element in $\mathfrak{p}_i$, for
each $i=1,\ldots, m$. Consider the polynomial
$$f(x):=(1+Mx)^3+aM^2(1+Mx)^2+bM^4(1+Mx)+cM^6~~\in A[x].$$

Then, since $k$ is a Hilbertian field by \cite[Ch.14,
Coro.~14.10]{jar}, $A$ is a Hilbertian ring. So there exists $m_1\in
A$ such that $f(m_1)$ is not a square in the completion $k_\infty$
of $k$ at $\infty$, so in particular, it is not a square in $k$. Let
$K_1=k(\sqrt{f(m_1)})$. Then, $\infty$ does not split in $K_1$ so
$K_1$ is an imaginary quadratic extension of $k$.

Next, by the Hilbert irreducibility over $K_1$ and by \cite[Ch.11,
Coro.~11.7]{jar}, we can choose $m_2\in A$ such that $f(m_2)$ is
non-square in both $K_1$ and $k_\infty$. Let $K_2=k(\sqrt{f(m_2)})$.
Then $K_1$ and $K_2$ are linearly disjoint imaginary quadratic
extensions over $k$.

 By repeating this procedure over the composite field
$K_{1}K_{1}\cdots K_{n}$ of imaginary quadratic extensions obtained
from the previous steps inductively, we obtain an infinite sequence
$\{K_i=k(\sqrt{f(m_i)})\}_{i=1}^{\infty}$ of quadratic extensions of
$k$ such that for all $i$,
\begin{enumerate}
\item $\infty$ does not split in $K_i$, and \item
 the fields
$K_{i}$ are pairwise linearly disjoint over $k$, $(
\mbox{\emph{i.e.} }[K_{1}K_{2}\cdots K_{r}:k]=2^r, \mbox{ for any }
r \geq 1 ).$
\end{enumerate}

Note that for every prime $\mathfrak{p}_j$ dividing $\mathfrak{n}$,
$$f(m_i) \equiv 1 ~~~ (\bmod~~\Fp_j).$$ So this implies that all
primes dividing $\mathfrak{n}$ split in $K_i$ for all $i$ (see
\cite[Proposition~10.5]{rosen}).

Let $\sigma \in G_k$. Then either $\sigma|_{K_i}=id_{K_i}$ for all
$i$, or $\sigma|_{K_i}\neq id_{K_i}$ for some $i$.

First, suppose that for all $i$, $\sigma|_{K_i}=id_{K_i}$. Then, for
each $i$, consider the element $\displaystyle \frac{1+Mm_i}{M^2}\in k$.
By plugging this
 into the given Weierstrass equation of $E/k$, we get $$
y^2=\left(\displaystyle\frac{1+Mm_i}{M^2}\right)^3
+a\left(\frac{1+Mm_i}{M^2}\right)^2+b\left(\frac{1+Mm_i}{M^2}\right)+c
=\frac{f(m_i)}{M^6}.$$ Hence, if we let
$$P_i=\left(\displaystyle\frac{1+Mm_i}{M^2},
\frac{\sqrt{f(m_i)}}{M^3}\right),$$ then  $P_i$ is a point in
$E(K_i)$ but it is not in $E(k)$.  And  moreover, since
$K_i=K_i^{\sigma}$, $P_i$ is fixed under $\sigma$.

So we get an infinite sequence $\{P_i\}_{i=1}^\infty$ of points in
$E(\ks)$ such that each $P_i$ is defined over the imaginary
quadratic extension $K_i$ over $k$. We may assume that these points
$P_i$ are not torsion points by Lemma~\ref{degreebd}. Now we show
the points $P_i$
 are linearly independent. Suppose that they are dependent. Then, for some
integers $a_j$ ,
$$(***) \hspace{5 cm} a_1P_{1}+a_2P_{2}+\cdots+a_rP_{r}=O.\hspace{5 cm}$$ Since the fields  $K_i$
  are pairwise linearly disjoint  over $k$, for each $i$, there
  is an automorphism of $\overline{k}$ which fixes all
  but one $K_i$ of $K_1,\ldots,K_r$. Note that such an automorphism takes
  $P_i$ to its inverse, $-P_i$. Applying this
  automorphism to $(***)$, we get
  $$a_1P_{1}+\cdots +a_{i-1}P_{i-1}-a_{i}P_{i}+\cdots +a_rP_{r}=O.$$
By subtracting this from $(***)$, we get $2a_iP_{i}=O$, which
implies $a_i=0$ since the characteristic $p$  of $k$ is not $2$ and
$P_i$ is not a torsion point. We conclude that the $P_i\in
E(\overline{k})$ are linearly independent. Moreover, $P_i$ are
defined over the composite field of all quadratic field extensions
of $k$, which is an abelian extension of $k$.
   Hence,  the rank of $E$ over
the maximal abelian extension of $k$ in $\ks$ is infinite,
   so the rank of $E(\ks)$  is infinite.

Next, suppose that there is an integer $i$ such that
$\sigma|_{K_i}\neq id _{K_i}$. Then, fix such a quadratic imaginary
extension $K_i$, and call it $K$. Then, our construction shows that
$K$ satisfies the hypothesis of Proposition~\ref{prop:sigma} (that
is, $(E,K)$ satisfies the Heegner hypothesis). So we complete the
proof of this case as a consequence of Proposition~\ref{prop:sigma}.
\end{proof}

\section{Modular elliptic curves over totally real number fields and infinite rank}

In this section we treat the case of elliptic curves parametrized by certain Shimura curves.

Let $F$ be a totally real number field with ring of integers $\CO_F$. For any abelian group $M$
we denote by
$\hat{M} = M\otimes\prod_p \bbz_p$ its profinite completion.
The ring of ad\`eles of $F$ is denoted by $\bba_F$.

Let $N\subset\CO_F$ be a non-zero ideal, and let $f$ be a newform on $\GL_2(\bba_F)$ of
parallel weight $2$, level
\[
K_0(N) = \left\{\big(\begin{smallmatrix}a & b \\ c &
d\end{smallmatrix}\big)\in\GL_2(\hat{\CO}_F) \;|\; c\in\hat{N}\right\},
\]
trivial central character, and rational Hecke eigenvalues.

Then, (see \cite{Zhang}) there exists an elliptic curve $E/F$ of
conductor $N$ such that
\begin{enumerate}
\item the $L$-functions of $E$ and $f$ coincide up to factors at primes dividing $N$, and
\item there exists a Shimura curve $X/F$ and a surjective $F$-morphism
\[
\pi : X \longrightarrow E.
\]
\end{enumerate}

We will refer to such elliptic curves as {\em modular} elliptic
curves, for example all elliptic curves over $F=\bbq$ are modular,
by the celebrated results of Wiles et al \cite{Wi,TW,BCDT}.


Our second main result is the following.


\begin{thm}\label{thm:main2}
Let $E/F$ be a modular elliptic curve of conductor $N$ over a
totally real number field $F$. If either $[F:\bbq]$ is odd, or $N$
is non-trivial, then for each $\sigma\in G_F$, the rank of
$E(\bar{F}^\sigma)$ is infinite.
\end{thm}

Our approach is similar to the function field case, and we first
show that $E$ has infinite rank over a suitable tower of ring class
fields.

Let $K/F$ be a totally imaginary quadratic extension, and denote by
\[
\varepsilon = \otimes_\nu\varepsilon_\nu :
F^\times\backslash\hat{F}^\times \longrightarrow \{\pm 1\}
\]
the character associated to $K/F$. We say $(E,K)$ satisfies the {\bf
weak Heegner Hypothesis} if
\begin{enumerate}
\item the relative discriminant of $K/F$ is prime to $N$, and
\item $\varepsilon(N) = (-1)^{[F:\bbq]-1}.$
\end{enumerate}

Throughout this section, we let $E/F$ be a modular elliptic curve
with conductor $N$ and $K/F$ a totally imaginary quadratic extension
such that $(E,K)$ satisfies the weak Heegner hypothesis. 

\begin{prop}\label{classfields}
Let $\Fp\subset\CO_F$ be a prime ideal not dividing $2N$ and satisfying $\varepsilon_\Fp(N)=1$,
and denote by $K[\Fp^n]$ the ring class field of $K$ with conductor $\Fp^n$.
Then the rank of $E(K[\Fp^n])$ is unbounded as $n\rightarrow\infty$.
\end{prop}

Let $K[\Fp^\infty]=\cup_{n=1}^{\infty}K[\Fp^n]$, and we start with
the following lemma.

\begin{lem}\label{finitetorsion}
$E(K[\Fp^\infty])$ has finite torsion.
\end{lem}

\begin{proof}
Let $\Fq_1$ and $\Fq_2$ be two principal primes of $F$ which are inert in $K/F$ and
at which $E$ has good reduction. Since they split completely in each $K[\Fp^n]$,
$K[\Fp^\infty]$ has finite residue fields $k_1$ and $k_2$ at $\Fq_1$ and $\Fq_2$,
respectively. From good reduction follows that $E(K[\Fp^{\infty}])_{\tor}$
injects into $E(k_1)\oplus E(k_2)$, which is finite.
\end{proof}

To prove Proposition \ref{classfields} we need to construct suitable Heegner points
on the Shimura curve $X$, for which we will need to introduce some notation.
Our standard reference is the article of Zhang \cite{Zhang}.

Fix a real place $\tau$ of $F$. Then there exists a unique quaternion algebra $B$
which is ramified precisely at $\tau$ and at all the finite places $\nu$
with $\varepsilon_\nu(N)=-1$ (the cardinality of this set of places is appropriate,
because we are assuming the weak Heegner hypothesis). We fix an embedding
\[
\rho : K \hookrightarrow B.
\]
Let $R\subset B$ be an order of type $(N,K)$, in other words $R$ contains
$\rho(\CO_K)$ and has conductor $N$. The Shimura curve $X/F$ corresponds to the Riemann surface
\[
X(\bbc) \cong B_+\backslash\bbh\times\hat{B}^\times/\hat{F}^\times \hat{R}^\times\cup\{\text{cusps}\},
\]
where $B_+$ denotes the elements of $B$ of totally positive reduced norm, $\bbh$
denotes the complex upper half-plane, and $\{\text{cusps}\}$ is a finite set,
which is non-empty only in the case where $F=\bbq$ and $X=X_0(N)$.

For the construction of Heegner points it is more convenient to work
with the Shimura curve $Y$ corresponding to the Riemann surface
\[
Y(\bbc) \cong B^\times\backslash\bbh^{\pm}\times\hat{B}^\times/\hat{R}^\times\cup\{\text{cusps}\},
\]
of which $X$ is a quotient by the action of $\hat{F}^\times$.

A point $z\in Y(\bbc)$ is called a {\em CM point} if it is represented by an element
of $\bbh^{\pm}\times\hat{B}$ of the form $(\sqrt{-1},g)$. To a CM point $z$ we associate the morphism
\[
\phi_z = g^{-1}\rho g : K \longrightarrow \hat{B}.
\]
The order $\End(z) := \phi_z^{-1}(\hat{R})$ in $K$ is called the {\em endomorphism ring} of $z$,
and does not depend on the choice of $g$. It is of the form
\[
\End(z) = \CO_F + c\,\CO_K,
\]
for an ideal $c\subset\CO_F$ called the {\em conductor} of $z$.

\begin{proof}[Proof of Proposition \ref{classfields}]
Let $\Fp\subset\CO_F$ be a prime ideal as in the statement of the proposition. Denote by $F_\Fp$
the completion of $F$ at $\Fp$ with 
uniformizer $\varpi$. $B$ splits at $\Fp$, and we choose an isomorphism
\[
B\otimes F_\Fp \cong M_2(F_\Fp)
\]
such that $\rho(\sqrt{-d})\otimes 1$ in $\rho(K)\otimes F_\Fp$ corresponds to
the matrix $\big(\begin{smallmatrix}0 & -1 \\ d & 0\end{smallmatrix}\big)\in M_2(F_\Fp)$,
where $K=F(\sqrt{-d})$, $d\in \CO_F$.

Now let $P\in\hat{B}^\times$ be the element with $\Fp$-component
$\big(\begin{smallmatrix}\varpi & 0 \\ 0 & 1\end{smallmatrix}\big)$ and all other components $1$.
Let $z_n$ be the CM point in $Y(\bbc)$ corresponding to
\[
(\sqrt{-1},P^n)\in\bbh^{\pm}\times\hat{B}^\times.
\]
As $\Fp\nmid 2N$ we see that $z_n$ has conductor $\Fp^n$, i.e.
\[
\End(z_n) = \CO_n := \CO_F + \Fp^n\CO_K.
\]

Denote by $x_n\in X(\bbc)$ and $y_n\in E(\bbc)$ the respective images of $z_n\in Y(\bbc)$
under the maps
\[
Y \longrightarrow X \stackrel{\pi}{\longrightarrow} E.
\]
We call the points $y_n$ {\em Heegner points} (in contrast, Zhang only uses the term Heegner
points for CM points with trivial conductor). Moreover, the points $x_n$, and thus also $y_n$,
are defined over $K[\Fp^n]$. In fact, by \cite[\S2.1.1]{Zhang} the set $X_n$ of (positively oriented)
CM points on $X$ with conductor $\Fp^n$ is in bijection with
$K^\times\backslash\hat{K}^\times/\hat{\CO}_n^\times\cong\Pic(\CO_n)$, with the action
by $\Gal(K[\Fp^n]/K)$ given by class field theory.

By Lemma \ref{finitetorsion} and \cite[Lemma~2.5]{im3}, it suffices
to show that the set
\[
\{y_n \;|\; n=1,2,\ldots\} \subset E(K[\Fp^\infty])
\]
is not finitely generated. Suppose it is finitely generated, then by the Mordell-Weil Theorem
it is contained in $E(L)$ for some number field $L$, which we may suppose contains $K$.
Let $d$ be the degree of the $F$-morphism $\pi : X \rightarrow E$. Then $G_L = \Gal(\bar{L}/L)$
acts on the fibers $\pi^{-1}(y_n)$, giving an upper bound for
the $G_L$-orbit of $x_n$: $\# G_L\cdot x_n \leq d$. On the other hand, we have the lower bound
$G_L\cdot x_n \geq \#\Pic(\CO_n)/[L:K]$, which is unbounded as $n\rightarrow\infty$.
This contradiction completes the proof of Proposition \ref{classfields}.
\end{proof}

Finally, the following proposition is analogous to
Proposition~\ref{prop:sigma}, hence completes our second main
result, Theorem~\ref{thm:main2}, whose proof is given below.

\begin{prop}\label{prop:sigma2} Let $\sigma\in G_F$.
Let $\Fp\subset\CO_F$ be a prime ideal not dividing $2N$ and
satisfying $\varepsilon_\Fp(N)=1$.
 Then, the rank of the Mordell-Weil group
$E((K[\Fp^n])^\sigma)$ over the fixed subfield of $K[\Fp^n]$ under
$\sigma$ is unbounded, as $n$ goes to $\infty$. In particular, the
rank of $E(K_{ab}^\sigma)$ is infinite, where $K_{ab}$ is the
maximal abelian extension of $K$.
\end{prop}

\begin{proof} By generalizing the result in \cite[Lemma 2.3]{im3} which is
elementary class field theory, we note that the Galois group of a
ring class field over $K$ of conductor $\Fp^n$ over $K$ is a product
of $p$-cyclic groups, where $p$ is the rational prime below $\Fp$,
and it has a dihedral group structure. So the proof is identical
with the argument in Proposition~\ref{prop:sigma} together with the
unboundedness of the rank of $E(K[\Fp^n])$ as $n$ goes to infinity
shown in Proposition~\ref{classfields}.
\end{proof}


We are now ready to prove Theorem \ref{thm:main2}. We need the
following two lemmas.

\begin{lem}\label{lem:dense} Let $K$ be a number field and
$\tau_1,\ldots,\tau_m$ be a family of real embeddings of $K$. For
$i=1,2,\ldots,k$, let $f_i(x,y)\in K[x,y]$ be irreducible
polynomials over $K(x)$. Let $H_K(f_i)$ $= \{\alpha\in K :$ $
f_i(\alpha, y) \in K[y]$ is irreducible over $K\}$ be the Hilbert
set of $f_i$ over $K$. Then  for any open interval $I$ in $\bbr$,
$$\left(\bigcap\limits _{i=1}^k
H_K(f_i)\right)~\cap~\left(\bigcap\limits
_{j=1}^m\tau_j^{-1}(I)\right)\neq~~ \emptyset.$$
\end{lem}

\begin{proof}
See  \cite[Lemma~1.2]{im3}.
\end{proof}

\begin{lem}\label{imag} Let $F$ be a totally real number field with
real embeddings $\tau_j$ for $j=1, \ldots, n$. Then for every
element $a\in F$ such that $\tau_j(a)<0$ for all $j$, the field
$F(\sqrt{a})$  is a totally imaginary quadratic extension over $F$.
\end{lem}

\begin{proof} Let $a\in F$ such that $\tau_j(a)<0$ for all
$j$. First, if $\sqrt{a}\in F$, then $\tau_j(\sqrt{a})\in\bbr$ for
all $j$ so its square $\tau_j(a)$ must be positive but this is not
true by the assumption. So $\sqrt{a}\notin F$ and the field
$F(\sqrt{a})$ is quadratic over $F$. Suppose $F(\sqrt{a})$ has a
real embedding $\rho$. Then, the restriction of $\rho$ to $F$ is one
of the real embeddings of $F$, so let $\rho|_F=\tau_k$ for some $k$.
Also, since $\rho(\sqrt{a})\in\bbr$,
$$0 < (\rho(\sqrt{a}))^2 = \rho(a)=\tau_k(a) <0,$$
which leads a contradiction. So $F(\sqrt{a})$ has no real embedding,
so it is totally imaginary.
\end{proof}

\begin{proof}[Proof of the main Theorem~\ref{thm:main2}]
Since $F$ has characteristic 0, we fix a Weierstrass equation of
$E/F$, $y^2=x^3+ax+b$. By a change of variables, we may assume
that $a$ and $b$ are algebraic integers in $F$. Since $F$ is a
totally real number field, there exist finitely many real embeddings
$\tau_j$ of $F$, for $j=1,\ldots, [F:\bbq]$. Also, we let $L$ be a
quadratic field extension of $F$ whose relative discriminant is $N$.

Let $M=p_1p_2\cdot\cdot\cdot p_k$, where $p_i$ are prime elements in
distinct prime ideals $\mathfrak{p_i}$ in $\mathcal{O}_F$ dividing
$N$. Consider the polynomial
$$f(x)=(1+Mx)^3+aM^4(1+Mx)+bM^6~~\in\mathcal{O}_F[x].$$
For a polynomial $g$ over $F$, denote by $g^{\tau_j}$ the polynomial
obtained by applying $\tau_j$ to each coefficient of $g$.

Then, there exists a real number $r$ such that for all $x<r$, the
expressions $f^{\tau_j}(x)$ for all $j$ are strictly negative. Let
$I=(-\infty, r)$ be the open interval in $\bbr$ of all real numbers
less than $r$. Choose an integer $m_1\in I\cap H_L(y^2-f(x))$. Then,
since $f^{\tau_j}(m_1)<0$ for all $j$ and $\sqrt{f(x)}$ is not in
$L$, $K_{m_1}:=\bbq(\sqrt{f(m_1)})$ is an imaginary quadratic
extension of $F$ by Lemma~\ref{imag} whose relative discriminant is
different from that of $L$.

   By Lemma~\ref{lem:dense}, there
exists a rational number in $I\cap (\bigcap\limits_{j=1}^n
H_{LK_{m_1}}(y^2-f(x)))$. In particular, since the Hilbert set
$H_{LK_{m_1}}(y^2-f(x))$ contains infinitely many rational primes by
\cite[Chapter~9, Corollary~2.4]{l83}, we can choose an integer
$m_2\in I\cap H_{LK_{m_1}}(y^2-f(x))$. Then,
$K_{m_2}:=\bbq(\sqrt{f(m_2)})$ is a quadratic imaginary extension of
$F$ by Lemma~\ref{imag}, and  $L$, $K_{m_1}$ and $K_{m_2}$ are
distinct, hence they are linearly disjoint over $F$. By repeating
this procedure over the composite field $LK_{m_1}K_{m_2}\cdots
K_{m_r}$ of imaginary quadratic extensions obtained from the
previous steps inductively, we obtain an infinite set $S$ of
integers  such that for all $m\in S$,
\begin{enumerate}
\item $f^{\tau_j}(m)<0$ for all $j$, so that $K_m:=\bbq(\sqrt{f(m)})$ is a totally imaginary quadratic
 extension of $F$, \item the fields in the infinite
sequence $\{K_{m}\}_{m\in S}$ are linearly disjoint over $F$,
$$( \mbox{\emph{i.e.} }[K_{m_1}K_{m_2}\cdots
K_{m_r}:\bbq]=2^r, \mbox{  for any } m_i\in S \mbox{ and  for every
} r\geq 1 )$$\item the relative discriminant of $K_m/F$ is different
from $N$ and \item $f(m) \equiv 1 ~(\mod \mathfrak{p}_i)$ for
distinct prime ideals $\mathfrak{p}_i$ in $\mathcal{O}_F$ dividing
$N$.
\end{enumerate}

Then for each $m\in S$, $(E,K_m)$ satisfies the first weak Heegner
hypothesis by (3) as above and it satisfies the second weak Heegner
hypothesis by (4) and by the assumption on $[F:\bbq]$ or $N$.

So for $\sigma\in G_F$, depending on whether $\sigma$ fixes $K_m$
for all $m\in S$ or not, the rest of the proof is the same as
\cite[Theorem~1.3]{im3} (or our first main result, Theorem
\ref{main}) by using Proposition~\ref{prop:sigma2}.
\end{proof}

\end{document}